\documentclass[12pt]{article}
\usepackage{amssymb}
\usepackage{latexsym,bm}
\usepackage{graphicx}
\usepackage{amsmath}
\usepackage{mathrsfs}
\usepackage{epstopdf}

\date{}
\newcounter{mathitem}
  {\begin{list}{{$(\roman{mathitem})$}}{
   \setcounter{mathitem}{0}
   \usecounter{mathitem}
   \setlength{\topsep}{0pt plus 2pt minus 0pt}
   \setlength{\parskip}{0pt plus 2pt minus 0pt}
   \setlength{\partopsep}{0pt plus 2pt minus 0pt}
   \setlength{\parsep}{0pt plus 2pt minus 0pt}
   \setlength{\leftmargin}{35pt}
   \setlength{\itemsep}{0pt plus 2pt minus 0pt}}}
  {\end{list}}

\begin{document}
\title{Sparse graphs with an independent or foresty minimum vertex cut
\footnote{E-mail addresses: {\tt chengkunmath@163.com} (K. Cheng), {\tt tyr2290@163.com} (Y. Tang), {\tt zhan@math.ecnu.edu.cn} (X. Zhan).}}
\author{\hskip -10mm Kun Cheng, Yurui Tang and Xingzhi Zhan\thanks{Corresponding author}\\
{\hskip -10mm \small Department of Mathematics, East China Normal University, Shanghai 200241, China}}\maketitle
\begin{abstract}
A connected graph is called fragile if it contains an independent vertex cut. In 2002 Chen and Yu proved that every connected graph of order $n$ and size
at most $2n-4$ is fragile, and in 2013 Le and Pfender characterized the non-fragile graphs of order $n$ and size $2n-3.$ It is natural to consider minimum vertex cuts.
We prove two results. (1) Every connected graph of order $n$ with $n\ge 7$ and size at most $\lfloor 3n/2\rfloor$ has an independent minimum vertex cut; (2) every connected graph of order $n$
with $n\ge 7$ and size at most $2n$ has a foresty minimum vertex cut. Both results are best possible.
\end{abstract}

{\bf Key words.} Fragile graph; minimum vertex cut; independent set; size; extremal problem

{\bf Mathematics Subject Classification.} 05C35, 05C40, 05C69
\vskip 8mm

\section{Introduction and main results}

We consider finite simple graphs and use standard terminology and notation from [1] and [9]. The {\it order} of a graph is its number of vertices, and the
{\it size} is its number of edges. We denote by $V(G)$ the vertex set of a graph $G,$ and for $S\subseteq V(G)$ we denote by $G[S]$ the subgraph of $G$ induced by
$S.$ A vertex cut of a connected graph $G$ is a set $S\subset V(G)$ such that $G-S$ is disconnected. A vertex cut $S$ of a connected graph $G$ is called an
{\it independent vertex cut} if $S$ is an independent set, and $S$ is called a {\it foresty vertex cut} if $G[S]$ is a forest. There is a recent work involving 
independent vertex cuts [6].

{\bf Definition 1} A connected graph is called {\it fragile} if it contains an independent vertex cut.

Fragile graphs have applications in some decomposition algorithms [2]. The following result was conjectured by Caro (see [4]) and proved by Chen and Yu [4] in 2002.

{\bf Theorem 1} [4] {\it Every connected graph of order $n$ and size at most $2n-4$ is fragile.}

The size bound $2n-4$ is sharp, and in 2013 Le and Pfender [7] characterized the non-fragile graphs of order $n$ and size $2n-3$ (see [8] for a related work).
Also in 2002  Chen, Faudree and Jacobson [3] proved the following result.

{\bf Theorem 2} [3] {\it Every connected graph of order $n$ and size at most $(12n/7)-3$ contains an independent vertex cut $S$ with $|S|\le 3.$}

Recently Chernyshev, Rauch and Rautenbach [5] has initiated the study of foresty vertex cuts of graphs. A vertex cut $S$ of a graph of connectivity $k$ is called
{\it minimum} if $|S|=k.$ It is natural to consider minimum vertex cuts.

In this paper we prove the following two results.

{\bf Theorem 3} {\it Every connected graph of order $n$ with $n\ge 7$ and size at most $\lfloor 3n/2\rfloor$ has an independent minimum vertex cut, and the size bound
$\lfloor 3n/2\rfloor$ is best possible.}

{\bf Theorem 4} {\it Every connected graph of order $n$ with $n\ge 7$ and size at most $2n$ has a foresty minimum vertex cut, and the size bound $2n$ is best possible.}

We give proofs of Theorems 3 and 4 in Section 2.

We denote by $|G|,$ $e(G)$ and $\kappa (G)$ the order, size and connectivity of a graph $G,$ respectively. The neighborhood of a vertex $x$ is denoted by $N(x)$ or $N_G(x),$
and the closed neighborhood of $x$ is $N[x]\triangleq N(x)\cup \{x\}.$ The degree of $x$ is denoted by ${\rm deg}(x).$ We denote by $\delta (G)$ and $\Delta(G)$ the minimum degree and maximum degree of $G,$ respectively.  For a vertex subset $S\subseteq V(G),$ we use $N(S)$ to denote the neighborhood  of $S;$ i.e., $N(S)=\{y\in V(G)\setminus S \,|\, y\,\,{\rm has}\,\,{\rm a}\,\,{\rm neighbor}\,\,{\rm in}\,\,S\}.$ For $x\in V(G)$ and $S\subseteq V(G),$ $N_S(x)\triangleq N(x)\cap S$ and the degree of $x$ in $S$ is ${\rm deg}_S(x)\triangleq |N_S(x)|.$ Given two disjoint vertex subsets $S$ and $T$ of $G,$ we denote by $[S, T]$ the set of edges having one endpoint in $S$ and the other in $T.$  The degree of $S$ is ${\rm deg}(S)\triangleq |[S, \overline{S}]|,$ where $\overline{S}=V(G)\setminus S.$ We denote by $C_n,$ $P_n$ and $K_n$ the cycle of order $n,$ the path of order $n$ and the complete graph of order $n,$ respectively. $\overline{G}$ denotes the complement of a graph $G.$ For two graphs $G$ and $H,$ $G\vee H$ denotes the {\it join} of $G$ and $H,$ which is obtained from the disjoint union $G+H$ by adding edges joining every vertex of $G$ to every vertex of $H.$

For graphs we will use equality up to isomorphism, so $G=H$ means that $G$ and $H$ are isomorphic.

\section{Proofs}

We will repeatedly use the following fact.

{\bf Lemma 5} {\it If $S$ is a minimum vertex cut of a connected graph $G,$ then every vertex in $S$ has a neighbor in every component of $G-S.$}

A $3$-regular graph is called a {\it cubic graph}.

{\bf Lemma 6} {\it Every connected cubic graph of order at least eight has an independent minimum vertex cut.}

{\bf Proof.} Let $G$ be a connected cubic graph of order at least $8.$ Then $\kappa(G)\in\{1,2,3\}.$ Lemma 6 holds trivially in the case $\kappa(G)=1.$  Next we
consider the remaining two cases.

Case 1. $\kappa(G)=2.$

Let $S=\{x, y\}$ be a minimum vertex cut of $G.$ If $x$ and $y$ are nonadjacent, then $S$ is what we want. Now suppose that $x$ and $y$ are adjacent. Let $H$ be a component
of $G-S.$ We assert that for any $v\in V(H),$ ${\rm deg}_S(v)\le 1.$ Otherwise $v$ would be a cut-vertex of $G,$ contradicting our assumption that $\kappa(G)=2.$
Since ${\rm deg}(x)=3$ and $x$ and $y$ are adjacent, $x$ has exactly one neighbor $p$ in $H.$ By the above assertion, $N_S(p)=\{x\},$ and consequently $p$ has two neighbors in $H.$
Then $\{p, y\}$ is an independent minimum vertex cut of $G.$

Case 2. $\kappa(G)=3.$

Choose a vertex $v\in V(G)$ and denote $S=N(v)=\{x,y,z\}.$ If $S$ is an independent set, then it is an independent minimum vertex cut of $G.$ Next suppose that
$S$ is not an independent set. Without loss of generality, suppose that $x$ and $y$ are adjacent. Since $G$ is cubic and $S$ is a minimum vertex cut, $\Delta(G[S])=1.$
It follows that $G[S]=K_2+K_1.$

Denote $T=V(G)\setminus N[v].$ We assert that for any $w\in T,$ $w$ is adjacent to at most one of $x$ and $y.$ Otherwise $\{w, z\}$ would be a vertex cut of $G,$ contradicting
our assumption $\kappa(G)=3.$ Let $\{p\}=N_T(x)$ and $\{q\}=N_T(y).$

We assert that at least one of $p$ and $q$ is nonadjacent to $z.$ To the contrary, suppose that both $p$ and $q$ are adjacent to $z.$ Since $G$ has order at least $8,$
$T\setminus \{p, q\}\neq\emptyset.$ Then $\{p,q\}$ is a vertex cut, contradicting our assumption $\kappa(G)=3.$ If $p$ is nonadjacent to $z,$ then $\{p, y, z\}$ is an independent minimum vertex cut of $G;$ if $q$ is nonadjacent to $z,$ then $\{q, x, z\}$ is an independent minimum vertex cut of $G.$ This completes the proof. \hfill $\Box$

The graph in Figure 1 shows that the lower bound $8$ for the order in Lemma 6 is sharp.
\begin{figure}[h]
\centering
\includegraphics[width=0.25\textwidth]{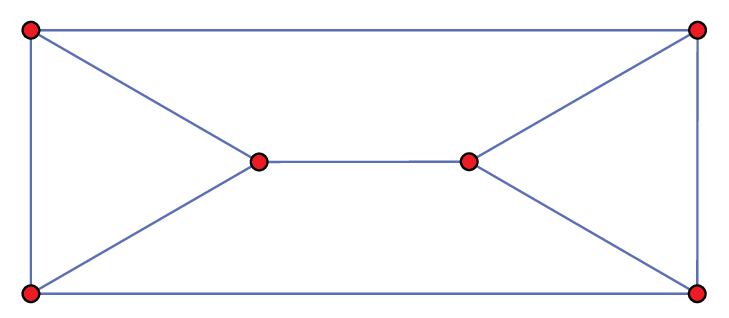}
\caption{A cubic graph of order $6$ without independent minimum vertex cut}
\end{figure}

{\bf Proof of Theorem 3.} We first use induction on the order $n$ to prove the statement that every connected graph of order $n$ with $n\ge 7$ and size at most $\lfloor 3n/2\rfloor$ has an independent minimum vertex cut.

{\bf The basis step.} $n=7.$

Let $F$ be a connected graph of order $7$ and size at most $10=\lfloor 3\times 7/2\rfloor.$ We have $\delta(F)\le 2,$ since otherwise we would have $e(F)\ge 11>10,$ a contradiction.
It follows that $\kappa(F)\le \delta(F)\le 2.$ The result holds trivially if $\kappa(F)=1.$ Thus it suffices to consider the case when $\kappa(F)=\delta(F)=2.$

Let $v$ be a vertex of degree $2$ and let $N(v)=\{x, y\}.$ If $x$ and $y$ are nonadjacent, then $\{x, y\}$ is an independent minimum vertex cut of $F.$ Next suppose that
$x$ and $y$ are adjacent. Applying Lemma 5 and using the size restriction of $F$ we deduce that $F-\{x, y\}$ has at most four components; i.e., $F-N[v]$ has at most three components. Then
$$
F-N[v]\in \{2K_1+K_2,\, 2K_2,\, K_1+P_3,\, K_1+C_3,\, K_1\vee\overline{K_3},\, P_4,\, C_4,\, K_1\vee (K_1+K_2),\, K_4^{-}\}
$$
where $K_4^{-}$ is the graph obtained from $K_4$ by deleting one edge. In each case it is easy to find two nonadjacent vertices which form a vertex cut by using the size condition.

{\bf The induction step.} $n\ge 8.$

Let $G$ be a connected graph of order $n\ge 8$ and size at most $\lfloor 3n/2\rfloor$ and suppose that the above statement holds for all graphs
of order $n-1.$ It suffices to consider the case $\kappa(G)\ge 2.$ Since $e(G)\le 3n/2,$ we have $2\le\kappa(G)\le\delta(G)\le 3.$

Case 1. $\delta(G)=3.$

Since $\delta(G)=3$ and $e(G)\le 3n/2,$ we have $\Delta(G)=3$ and hence $G$ is cubic. The statement holds by Lemma 6.

Case 2. $\delta(G)=2.$

In this case $\kappa(G)=2.$ Choose a vertex $v$ of degree $2$ and let $N(v)=\{x, y\}.$ If $x$ and $y$ are nonadjacent, then $\{x, y\}$ is an independent minimum vertex cut.
Next we assume that $x$ and $y$ are adjacent. Denote $H=G-v.$ Then $H$ is a connected graph of order $n-1$ and
$$
e(H)=e(G)-2\le \frac{3n}{2}-2=\frac{3n-4}{2}<\frac{3(n-1)}{2},
$$
which implies that $\delta(H)\le 2$ and hence $\kappa(H)\le 2.$ On the other hand, since $x$ and $y$ are adjacent, the condition $\kappa(G)=2$ implies that $\kappa(H)\ge 2.$
Thus $\kappa(H)=2.$ By the induction hypothesis, $H$ has an independent vertex cut $M$ with $|M|=2.$ Clearly $M$ is an independent minimum vertex cut
of $G.$

Now for every integer $n\ge 7$ we construct a graph $G_n$ of order $n$ and size $\lfloor 3n/2\rfloor+1$ such that $G_n$ has no independent minimum vertex cut. Hence
the size bound $\lfloor 3n/2\rfloor$ in Theorem 3 is best possible.

If $n$ is odd, let $C: v_1v_2\ldots v_{n-1}v_1$ be an $(n-1)$-cycle. Add a vertex $v_n$ to $C$ and then add edges $v_1v_{(n+1)/2},$ $v_1v_n,$ $v_{(n+1)/2}v_n,$
$v_iv_{n+1-i}$ for $i=2,3,\ldots,(n-1)/2$ to obtain $G_n.$ If $n$ is even, let $D: v_1v_2\ldots v_nv_1$ be an $n$-cycle. Then in $D$ add edges
$v_2v_{n/2},$ $v_{(n+4)/2}v_n,$ $v_iv_{n+2-i}$ for $i=2,3,\ldots,n/2$ to obtain $G_n.$ We depict $G_{11}$ and $G_{12}$ in Figure 2.
\begin{figure}[h]
\centering
\includegraphics[width=0.85\textwidth]{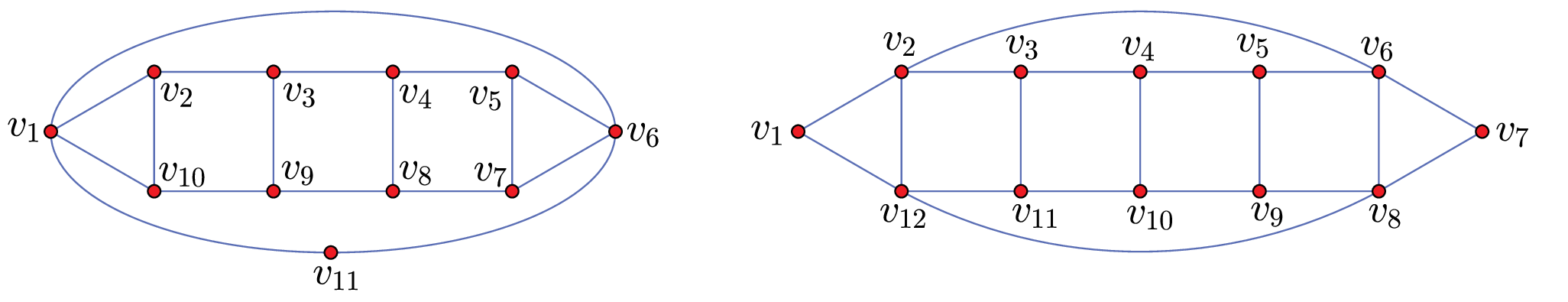}
\caption{$G_{11}$ and $G_{12}$}
\end{figure}

$G_n$ has order $n$ and size $\lfloor 3n/2\rfloor+1.$ If $n$ is odd, $\{v_1, v_{(n+1)/2}\}$ is the unique minimum vertex cut of $G_n,$ which induces an edge.
If $n$ is even, $G_n$ has exactly two minimum vertex cuts: $\{v_2, v_n\}$ and $\{v_{n/2}, v_{(n+4)/2}\},$ each of which induces an edge. Thus $G_n$ has no independent
minimum vertex cut. This completes the proof. \hfill $\Box$

Now we prepare to prove Theorem 4.

Let $S$ and $T$ be two disjoint vertex subsets of a graph $G.$ An {\it $(S, T)$-path} is a path $P$ with one endpoint in $S$ and the other in $T$ such that
$S\cup T$ contains no internal vertex of $P.$ The following fact is well-known [9, p.174] and it follows from Menger's theorem ([1, p.208] or [9, p.167]).

{\bf Lemma 7} {\it Let $G$ be a $k$-connected graph. If $S$ and $T$ are two disjoint subsets of $V(G)$ with cardinality at least $k,$ then $G$ has $k$ pairwise vertex disjoint
$(S, T)$-paths.}

A {\it $k$-matching} is a matching of cardinality $k.$

{\bf Lemma 8} {\it Let $S$ be a vertex cut of a $k$-connected graph $G$ and let $H$ be a component of $G-S.$ If $|H|\ge k,$ then the set $[S,\, V(H)]$
contains a $k$-matching.}

{\bf Proof.} Since $G$ is $k$-connected, $|S|\ge k.$ By Lemma 7, $G$ contains $k$ pairwise vertex disjoint $(S, V(H))$-paths $P_i,$ $i=1,\ldots,k.$
Clearly each $P_i$ must be an edge, and hence $\{P_1,P_2,\ldots,P_k\}$ is a $k$-matching in $[S,\, V(H)].$ \hfill $\Box$

{\bf Lemma 9} {\it Every connected $4$-regular graph of order at least seven has a foresty minimum vertex cut.}

{\bf Proof.} Let $G$ be a $4$-regular graph of order $n$ with $n\ge 7.$ We will show that $G$ has a foresty minimum vertex cut.
We have $\kappa(G)\le 4.$ If $\kappa(G)\le 2,$ the result holds trivially. Next suppose $\kappa(G)\ge 3$ and we distinguish two cases.

Case 1. $\kappa(G)=3.$

Let $S$ be a vertex cut of $G$ with $|S|=3.$ If $G[S]\neq C_3,$ then $S$ is a foresty minimum vertex cut of $G.$ Suppose $G[S]=C_3.$ Since $G$ is $4$-regular,
by Lemma 5 we deduce that $G-S$ has exactly two components, which we denote by $G_1$ and $G_2.$ Without loss of generality, suppose $|G_1|\ge |G_2|.$ Then
$|G_1|\ge (n-|S|)/2\ge (7-3)/2=2.$ Let $S=\{x,y,z\}.$ We assert that $d_S(v)\le 1$ for any $v\in V(G_1).$ To the contrary, suppose that there is $v\in V(G_1)$ such that
$d_S(v)\ge 2.$ Without loss of generality, suppose $\{x, y\}\subseteq N_S(v).$ Then $\{v, z\}$ is a vertex cut of $G,$ contradicting the assumption that  $\kappa(G)=3.$

Let $u$ be the neighbor of $x$ in $G_1.$ Then $\{u, y, z\}$ is a vertex cut of $G$ which induces $K_1+K_2.$ Hence it is a foresty minimum vertex cut.

Case 2. $\kappa(G)=4.$

Choose a vertex $v\in V(G)$ and denote $T=N(v)=\{x,y,z,u\}.$ Then $T$ is a minimum vertex cut of $G.$ Denote $H=G[T].$  If $H$ is a forest, then $T$ is what we want. Next suppose $H$ contains a cycle. By Lemma 5 and the condition that $G$ is $4$-regular, we have $\Delta(H)\le 2.$ Thus $H\in \{C_4, C_3+K_1\}.$

Subcase 2.1. $H=C_4.$

Let $W=V(G)\setminus N[v].$ We assert that for any $w\in W,$ $d_T(w)\le 1.$ Otherwise there exists a $w\in W$ with $d_T(w)\ge 2.$ Since the order $n\ge 7$ and $G$ is $4$-regular, $N(w)\neq T.$ Now $\{w\}\cup T\setminus N(w)$ is a vertex cut of cardinality at most $3,$ contradicting $\kappa(G)=4.$

Since $G$ is $4$-regular, by Lemma 5 we deduce that every vertex in $T$ has exactly one neighbor in $W.$ Let $f$ be the neighbor of $x$ in $W.$ Then
$R\triangleq \{f,y,z,u\}$ is a minimum vertex cut of $G$ and $G[R]=K_1+P_3$ is a forest.

Subcase 2.2. $H=C_3+K_1.$

Without loss of generality, suppose that $G[A]=C_3$ where $A=\{x,y,z\}.$ We assert that every vertex in $W$ has at most one neighbor in $A.$ Otherwise, there exists a vertex
$w\in W$ which has at least two neighbors in $A.$ Then $\{w, u\}\cup A\setminus N(w)$ is a vertex cut of $G$ of cardinality at most $3,$ contradicting $\kappa(G)=4.$

Let $p$ be the neighbor of $x$ in $W.$ Then $\{p, y, z, u\}$ is a foresty minimum vertex cut of $G.$ \hfill $\Box$

{\bf Remark.} There is only one $4$-regular graph of order $6,$ which has connectivity $4$ and has no foresty minimum vertex cut. Thus the lower bound $7$ for the order
in Lemma 9 is sharp.

{\bf Proof of Theorem 4.}  The first part of Theorem 4 is the following

{\bf Statement.} Every connected graph of order $n$ with $n\ge 7$ and size at most $2n$ has a foresty minimum vertex cut.

We use induction on the order $n$ to prove this statement. 

{\bf The basis step.} $n=7.$

Let $M$ be a graph of order $7$ and size at most $14.$ The condition $e(M)\le 14$ implies $\kappa(M)\le \delta(M)\le 4.$ If $\kappa(M)\le 2$ then the statement holds trivially.
 Next suppose $3\le\kappa(M)\le\delta(M)\le 4.$

If $\delta(M)=4,$ then $M$ is $4$-regular and by Lemma 9, $M$  has a foresty minimum vertex cut. Now suppose $\delta(M)=3.$ Then $\kappa(M)=3.$ Choose a vertex $v\in V(M)$ with 
${\rm deg}(v)=3,$  let $S=N(v)$ and let $R=V(M)\setminus N[v].$ If $M[S]$ is a forest, then $S$ is a foresty minimum vertex cut. Now suppose that $M[S]=C_3.$
If $R$ is an independent set, then the condition $\delta(M)=3$ implies that $e(M)=15,$ contradicting $e(M)\le 14.$ Hence $M[R]\in \{K_2+K_1, P_3, C_3\}.$ In each of these three cases,
 it is easy to verify that the statement holds by applying Lemma 8 in the latter two cases.

{\bf The induction step.} $n\ge 8.$

Let $G$ be a connected graph of order $n$ with $n\ge 8$ and size at most $2n,$ and suppose that the above statement holds for all graphs of order $n-1.$
The condition $e(G)\le 2n$ implies $\kappa(G)\le \delta(G)\le 4.$ If $\kappa(G)\le 2$ then the statement holds trivially.
Next suppose $3\le\kappa(G)\le\delta(G)\le 4.$

Case 1. $\delta(G)=4.$

Since $e(G)\le 2n,$ $G$ is $4$-regular. The statement holds by Lemma 9.

Case 2. $\delta(G)=3.$

We have $\kappa(G)=3.$ Choose a vertex $v\in V(G)$ with ${\rm deg}(v)=3$ and denote $S=N(v).$ If $G[S]$ is a forest, then $S$ is a foresty minimum vertex cut. Otherwise
$G[S]=C_3.$ Consider the graph $H=G-v.$ $H$ has order $n-1$ and  $e(H)=e(G)-3\le 2n-3<2(n-1),$ which implies that $\delta(H)\le 3.$ Hence $\kappa(H)\le 3.$ Since any vertex cut
of $H$ is also a vertex cut of $G$ and $\kappa(G)=3,$ we deduce that $\kappa(H)=3.$ By the induction hypothesis, $H$ has a foresty minimum vertex cut $T.$
Clearly $T$ is also a foresty minimum vertex cut of $G.$

\begin{figure}[h]
\centering
\includegraphics[width=0.8\textwidth]{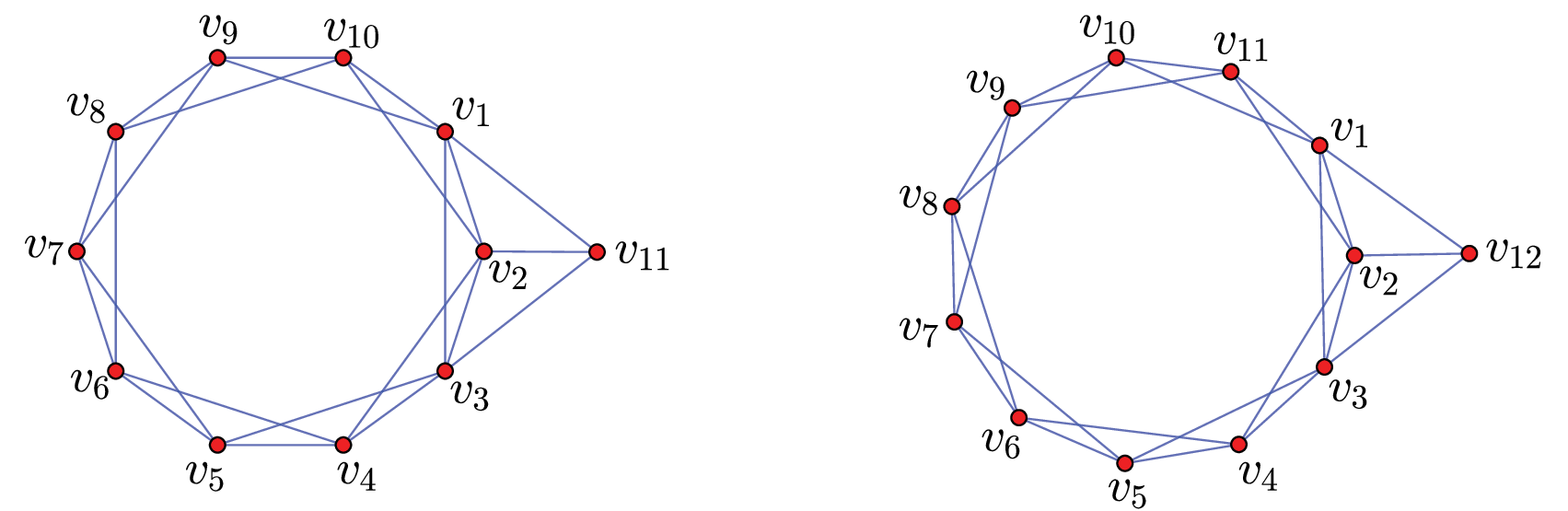}
\caption{$F_{11}$ and $F_{12}$}
\end{figure}

Now for every integer $n\ge 7$ we construct a graph $F_n$ of order $n$ and size $2n+1$ such that $F_n$ has no foresty minimum vertex cut. This shows that the size bound $2n$
in Theorem 4 is best possible. Recall that a chord $xy$ of a cycle $D$ is called a $k$-chord if the distance between $x$ and $y$ on $D$ is $k.$ Let $C: v_1v_2\ldots v_{n-1}v_1$
be a cycle of order $n-1.$ Add all the $2$-chords to $C$ to obtain a $4$-regular graph $R.$ Finally adding a new vertex $v_n$ to $R$ and adding the edges $v_nv_1,$ $v_nv_2$ and
$v_nv_3,$ we obtain $F_n.$ We depict $F_{11}$ and  $F_{12}$ in Figure 3.

It is easy to see that $\kappa(F_n)=3$ and $\{v_1, v_2, v_3\}$ is the unique minimum vertex cut, which induces a triangle. \hfill $\Box$

{\bf Data availability}

No data was used for the research described in the article.

\vskip 5mm
{\bf Acknowledgement.} This research  was supported by the NSFC grant 12271170 and Science and Technology Commission of Shanghai Municipality
 grant 22DZ2229014.

\end{document}